\documentclass[11pt,epsfig]{amsart}
\usepackage{amsmath, amscd, amssymb}
\usepackage{graphicx, psfrag}
\usepackage[frame,cmtip,arrow,matrix,line,graph,curve]{xy}
\usepackage{latexsym, amssymb, amsmath, epsfig, amscd, amsfonts}
\usepackage{epsfig}
\usepackage{graphpap, color}
\usepackage[mathscr]{eucal}
\usepackage{mathrsfs}
\usepackage{pstricks}
\usepackage{color}
\usepackage{cancel}
\usepackage{concrete,euler,textcomp}
\usepackage{youngtab}
\usepackage[T1]{fontenc}

\numberwithin{equation}{section}
\newtheorem{theorem}{Theorem}[section]

\newtheorem{corollary}[theorem]{Corollary}
\newtheorem{definition}[theorem]{Definition}
\newtheorem{example}[theorem]{Example}
\newtheorem{lemma}[theorem]{Lemma}
\newtheorem{proposition}[theorem]{Proposition}
\newtheorem{remark}[theorem]{Remark}

\begin{document}
\title[Nullcone for the symplectic group]{The nullcone in the multi-vector
representation of the symplectic group and related combinatorics}
\author{Sangjib Kim}

\thanks{ The author is supported by iBAEKSU 0102-30-81.}

\address{ Department of Mathematics\\
The University of Arizona\\
617 N. Santa Rita Ave. P.O. Box 210089\\
Tucson, AZ 85721 USA}
\email{sangjib@math.arizona.edu}

\begin{abstract}
We study the nullcone in the multi-vector representation of the symplectic
group with respect to a joint action of the general linear group and the
symplectic group. By extracting an algebra over a distributive lattice
structure from the coordinate ring of the nullcone, we describe a toric
degeneration and standard monomial theory of the nullcone in terms of double
tableaux and integral points in a convex polyhedral cone.
\end{abstract}

\subjclass[2000]{20G05, 05E15, 13P10}
\keywords{Symplectic group, Nullcone, Young tableaux, Gelfand-Tsetlin
patterns, Toric deformations}
\maketitle

\section{Introduction}

Let $GL_{n}$ and $Sp_{2n}$ denote respectively the general linear group and
the symplectic group over the complex number field $\mathbb{C}$.
Combinatorics of tableaux provides a unifying scheme to understand
representation theory of $GL_{n}$ and geometry of the flag varieties and the
Grassmann varieties. In particular, the theory of double tableaux gives a
finite presentation of the coordinate ring of the affine space $M_{n,m}\cong 
\mathbb{C}^{n}\otimes \mathbb{C}^{m}$ which is compatible with the natural
action of $GL_{n}\times GL_{m}$. Moreover, we can explicitly describe weight
bases of representation spaces from the combinatorial structure of tableaux.

In this paper, we develop a parallel theory for the $Sp_{2n}$-nullcone $%
\mathcal{N}_{k,2n}$ which is defined by $Sp_{2n}$-invariant polynomials on
the space $M_{k,2n}$ with vanishing constant terms. We begin with a known
algebro-combinatorial description of the space $M_{n,m}$ as a cell of the
Grassmann variety of $n$ dimensional spaces in $\mathbb{C}^{m+n}$. Using
this observation, we construct a convex polyhedral cone $\mathcal{C}%
(M_{n,m}) $ associated with the space $M_{n,m}$ and study the integral
points in the cone. Then we characterize the defining ideal of $\mathcal{N}%
_{k,2n}$ in terms of integral points in $\mathcal{C}(M_{k,2n})$. This
characterization provides a convex polyhedral cone $\mathcal{C}(\mathcal{N}%
_{k,2n})$ associated with $\mathcal{N}_{k,2n}$ compatible with the action of 
$GL_{k}\times Sp_{2n}$. Our construction of the polyhedral cone $\mathcal{C}(%
\mathcal{N}_{k,2n})$ turns out to be related to a fiber product of the
Gelfand-Tsetlin patterns.

We also describe explicit joint weight vectors of $GL_{k}\times Sp_{2n}$ in
the coordinate ring of $\mathcal{N}_{k,2n}$ in terms of standard double
tableaux. As a result, we obtain standard monomial theory for the nullcone
and show that the nullcone can be degenerated to an affine toric variety
presented by an algebra over a distributive lattice.

\medskip

The toric degenerations of spherical varieties (e.g., \cite{AB04, Ca02,
GL96, Stu95}) and standard monomial theory (e.g., \cite{La03, Mu03}) have
been actively studied. Using classical invariant theory, we can study such
combinatorial and geometric results in connection with various
representation theoretic questions.

The recent papers \cite{HJLTW, HL1, HTW2} and their sequels construct
algebras encoding branching rules of representations of the classical
groups, and then study their standard monomial bases and toric
degenerations. With a similar philosophy, \cite{HTW3} and \cite{KL} study
tensor products of representations for the classical groups with explicit
highest weight vectors. By degenerating the muti-homogeneous coordinate
rings of the flag varieties, \cite{Ki08} and \cite{Ki2} describe weight
vectors of the classical groups in terms of the Gelfand-Tsetlin polyhedral
cone. This paper may be understood as an application of such approaches to
the nullcone in the multi-vector representation of the symplectic group to
obtain explicit combinatorial and representation theoretic descriptions of
it. For the nullcone associated with representations of reductive groups and
its interesting applications, we refer readers to \cite{KW06, KW2}.

\medskip

This paper is arranged as follows: In Section \ref{Mnm}, we introduce
notations for tableaux and patterns, and review standard monomial theory for
the coordinate ring $\mathbb{C}[M_{n,m}]$ of $M_{n,m}$. In Section \ref{CMnm}%
, we define the convex polyhedral cone $\mathcal{C}(M_{n,m})$ associated
with a degeneration of $M_{n,m}$, and study its connection to
representations of the general linear group. In Section \ref{Nk2n}, we study
the coordinate ring $\mathcal{R}(\mathcal{N}_{k,2n})$ of $\mathcal{N}_{k,2n}$
and show its standard monomial theory and toric degeneration. In Section \ref%
{CNk2n}, we describe the integral points in the convex polyhedral cone $%
\mathcal{C}(\mathcal{N}_{k,2n})$ for $\mathcal{N}_{k,2n}$ and explain its
relations to representation theory.

\medskip

\section{Affine space: $M_{n,m}$\label{Mnm}}

In this section, we review some results on the Young tableaux, the
Gelfand-Tsetlin patterns, and their applications to representation theory
and geometry of the Grassmann varieties.

\subsection{Tableaux}

Let $M_{n,m}=M_{n,m}(\mathbb{C})$ be the space of complex $n$ by $m$
matrices:%
\begin{equation}
M_{n,m}=\{(x_{ij}):1\leq i\leq n,1\leq j\leq m\}.  \label{coordinates}
\end{equation}

A \textit{Young diagram} or \textit{shape} is an array of square boxes
arranged in left-justified horizontal rows with each row no longer than the
one above it (e.g., \cite{FH91, Sta99}). We identify a shape $D$ with its
sequence of row lengths $D=(r_{1},r_{2},\cdots )$. Then the \textit{transpose%
} $D^{t}$ of $D$ is a shape $(c_{1},c_{2},\cdots )$ where $c_{i}$ is the
length of the $i$-th column of $D$. The \textit{length} $\ell (D)$ of shape $%
D$ is the number of rows in $D$. For example, the shape $%
(4,2,1)=(3,2,1,1)^{t}$ of length $3$ can be drawn as%
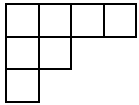
\begin{equation*}
\young(\ \ \ \ ,\ \ ,\ )
\end{equation*}

For a subset $I=[i_{1},\cdots ,i_{l}]$ (respectively $J=[j_{1},\cdots
,j_{l}] $) of $\{1,\cdots ,n\}$ (respectively $\{1,\cdots ,m\}$), which we
can think of a filling of shape $(1,\cdots ,1)$ of length $l$ with its
elements, we shall call the pair $[I:J]$ an\textit{\ one-line tableau} of
length $\ell ([I:J])=l$. We assume that the entries of $I$ and $J$ are
listed in increasing order, i.e., $1\leq i_{1}<\cdots <i_{l}\leq n$ and $%
1\leq j_{1}<\cdots <j_{l}\leq m$.

A partial order $\preceq $, called the \textit{tableau order}, can be
imposed on the set of one-line tableaux%
\begin{equation*}
D(n,m)=\{[I:J]:|I|=|J|\leq \min (n,m)\}
\end{equation*}%
as follows: $[I:J]\preceq \lbrack I^{\prime }:J^{\prime }]$, if the length
of $[I:J]$ is not smaller than that of $[I^{\prime }:J^{\prime }]$, and $%
i_{k}\leq i_{k}^{\prime }$ and $j_{k}\leq j_{k}^{\prime }$ for each $k$ not
bigger than the length of $[I^{\prime }:J^{\prime }]$. Then it is easy to
see that with respect to the tableau order $D(n,m)$ forms a distributive
lattice $(D(n,m),\wedge ,\vee )$.

Consider a collection $\{[I_{1}:J_{1}],\cdots ,[I_{u}:J_{u}]\}\subset D(n,m)$
with $l_{k}=\ell ([I_{k}:J_{k}])$ for each $k$. A concatenation $\mathsf{t}$
of its elements is called a \textit{double tableau, }if they are arranged so
that $l_{k}\geq l_{k+1}$ for all $k$. The \textit{shape} $sh(\mathsf{t})$ of 
$\mathsf{t}$ is the Young diagram $(l_{1},\cdots ,l_{u})^{t}$. 
For example,
the following one-line tableaux:%
\begin{equation*}
\lbrack 123:124],[13:23],[2:3],[4:5].
\end{equation*}
form a double tableau of shape $(4,2,1)$. 

We note that by considering first components $I_{k}$ and the second
components $J_{k}$ separately, we can think of a double tableau of shape $D$
in terms of a pair of fillings of the same shape $D$. 
The double tableau consisting of $[123:124]$, $[13:23]$, $[2:3]$ and $[4:5]$
can be matched with the following fillings $T^{-}$ and $T^{+}$ of shape $%
(4,2,1)$:%
\begin{equation}
\young(1124,23,3)\text{ and }\young(1235,23,4)  \label{splitting}
\end{equation}

Let us write $\delta _{\lbrack I:J]}$ for the map from $M_{n,m}$ to $\mathbb{%
C}$ by assigning to a matrix $X\in M_{n,m}$ the determinant of the $l\times
l $ submatrix of $X$ with rows and columns indexed by $I$ and $J$
respectively:%
\begin{equation*}
\delta _{\lbrack I:J]}=\det \left[ 
\begin{array}{ccc}
x_{i_{1}j_{1}} & \cdots & x_{i_{1}j_{l}} \\ 
\vdots & \ddots & \vdots \\ 
x_{i_{l}j_{1}} & \cdots & x_{i_{l}j_{l}}%
\end{array}%
\right]
\end{equation*}%
For a double tableau $\mathsf{t}$ consisting of $\left\{
[I_{k}:J_{k}]\right\} $, we define the corresponding element in the
coordinate ring $\mathbb{C}[M_{n,m}]$ of $M_{n,m}$ to be the following
product:%
\begin{equation*}
\Delta (\mathsf{t})=\prod_{1\leq k\leq u}\delta _{\lbrack I_{k}:J_{k}]}
\end{equation*}

\begin{definition}
\label{standardmonomials}A double tableau $\mathsf{t}$ is called a \textit{%
standard tableau} if the one-line tableaux in $\mathsf{t}$ form a multiple
chain in $D(n,m)$, i.e., 
\begin{equation*}
\lbrack I_{1}:J_{1}]\preceq \cdots \preceq \lbrack I_{u}:J_{u}]
\end{equation*}%
For a standard tableau $\mathsf{t}$, we call $\Delta (\mathsf{t})$ a \textit{%
standard monomial}.
\end{definition}

\medskip

Now, let us consider the following subposet:%
\begin{eqnarray*}
L(n,m) &=&\{[I:J]\in D(n,m):I=[1,\cdots ,|J|]\} \\
&\cong &\{[j_{1},\cdots ,j_{l}]:l\leq \min (n,m),1\leq j_{1}<\cdots
<j_{l}\leq m\}
\end{eqnarray*}%
For each $[I:J]\in L(n,m)$, since the first component $I$ is determined by
the size of $J$, we shall write $J$ for $[I:J]$. Similarly, we can also
think of the subposet $L^{\prime }(n,m)$ consisting of $[I:J]$ with $%
J=[1,\cdots ,|I|]$. Then, a \textit{semistandard tableau} as found in the
literature (e.g., \cite{FH91, Sta99}) can be defined as a multiple chain in
such posets with respect to the tableau order. For example, the
tableaux in (\ref{splitting}) are semistandard and they can be understood as
chains%
\begin{eqnarray*}
\lbrack 1,2,3] \preceq [1,3]\preceq \lbrack 2]\preceq \lbrack 4] \\
\lbrack 1,2,4] \preceq [2,3]\preceq \lbrack 3]\preceq \lbrack 5]
\end{eqnarray*}%
in $L^{\prime }(n,m)$ and $L(n,m)$ respectively.

This poset $L(n,m)$ has been extensively studied for the flag varieties, the
Grassmann varieties, and the determinantal varieties. For example, the
elements of $L(n,m+n)$ with fixed length $n$, which we shall denote by $%
Pl(n,m+n)$, may encode the Pl\"{u}cker coordinates for the Grassmann variety 
$Gr(n,m+n)$ of $n$ dimensional spaces in $\mathbb{C}^{m+n}$. In fact, this
case is general enough to study double tableaux thanks to the following
correspondence.

\begin{lemma}
\label{double2single}The following map $\xi $ gives an order isomorphism
from $D(n,m)$ to the subposet of $L(n,m+n)$ consisting of all the elements
of length $n$ except $[m+1,\cdots ,m+n]$: for $[I:J]\in D(n,m)$ with $%
I=[i_{1},i_{2},\cdots ,i_{h}]$ and $J=[j_{1},j_{2},\cdots ,j_{h}]$,%
\begin{equation*}
\xi :[I:J]\longmapsto \lbrack j_{1},j_{2},\cdots
,j_{h},m+u_{1},m+u_{2},\cdots ,m+u_{n-h}]
\end{equation*}%
where $\left\{ u_{k}\right\} $ is defined so that $\{n+1-u_{1},\cdots
n+1-u_{h}\}$ is complement to $I$ in $\{1,2,\cdots ,n\}$.
\end{lemma}

For the proof, see \cite[p. 519]{Pr07}.

\subsection{Gelfand-Tsetlin Patterns}

The poset $\widehat{L}(m,m)$, $L(m,m)$ with an extra top element, with
respect to the tableau order $\preceq $ turns out to be a distributive
lattice whose join-irreducibles form the following poset, which we shall
call the \textit{Gelfand-Tsetlin(GT) poset}:%
\begin{equation*}
\Gamma _{m}=\{z_{j}^{(i)}:1\leq j\leq i\leq m\}
\end{equation*}%
satisfying $z_{j}^{(i+1)}\geq z_{j}^{(i)}\geq z_{j+1}^{(i+1)}$ for all $i$
and $j$\ (\cite[Theorem 3.8]{Ki08}). We call $%
z^{(i)}=(z_{1}^{(i)},z_{2}^{(i)},\cdots ,z_{i}^{(i)})$ the $i$\textbf{-}%
\textit{th row} of $\Gamma _{m}$. We will draw it in a reversed triangular
array so that $z_{j}^{(i)}$ are decreasing from left to right along
diagonals. For example, $\Gamma _{4}$ can be drawn as%
\begin{equation*}
\begin{array}{ccccccc}
z_{1}^{(4)} &  & z_{2}^{(4)} &  & z_{3}^{(4)} &  & z_{4}^{(4)} \\ 
& z_{1}^{(3)} &  & z_{2}^{(3)} &  & z_{3}^{(3)} &  \\ 
&  & z_{1}^{(2)} &  & z_{2}^{(2)} &  &  \\ 
&  &  & z_{1}^{(1)} &  &  & 
\end{array}%
\end{equation*}

\begin{definition}
A GT pattern $\mathsf{p}$ of $GL_{m}$ is an order preserving map from $%
\Gamma _{m}$ to the set of non-negative integers:%
\begin{equation*}
\mathsf{p}:\Gamma _{m}\rightarrow \mathbb{Z}_{\geq 0},
\end{equation*}%
and the $i$-th row of $\mathsf{p}$ is $\mathsf{p}(z^{(i)})=(\mathsf{p}%
(z_{1}^{(i)}),\mathsf{p}(z_{2}^{(i)}),\cdots ,\mathsf{p}(z_{i}^{(i)}))$.
\end{definition}

The $m$-th row of a GT pattern $\mathsf{p}$ of $GL_{m}$ will be
alternatively called the \textit{top row} of $\mathsf{p}$. By identifying $%
\mathsf{p}$\ with its values, our definition agrees with the original one in 
\cite{GT50}.

Note that for each $i$, since $\mathsf{p}(z_{1}^{(i)})\geq \mathsf{p}%
(z_{2}^{(i)})\geq \cdots \geq \mathsf{p}(z_{i}^{(i)})\geq 0$, the $i$-th row 
$\mathsf{p}(z^{(i)})$ of a pattern $\mathsf{p}$ can be seen as a Young
diagram. There is a well-known conversion procedure between semistandard
Young tableaux and GT patterns compatible with successive branching rules of 
$GL_{k}$ down to $GL_{k-1}$ for $2\leq k\leq m$. See, e.g., \cite[%
Proposition 8.1.6]{GW98}.

\begin{lemma}
\label{pattern2tableau}The following procedure gives a bijection from the
set of GT patterns of $GL_{m}$ with a fixed $m$-th row $D$ and the set of
semistandard tableaux of shape $D$ with entries from $\{1,\cdots ,m\}$: for
a given GT pattern $\mathsf{p}$ of $GL_{m}$ with $\mathsf{p}%
(z^{(m)})=(r_{1},\cdots ,r_{m})$, fill in the cells in Young diagram $%
(r_{1},\cdots ,r_{m})$ corresponding to the skew diagram $\mathsf{p}%
(z^{(i)})/\mathsf{p}(z^{(i-1)})$ with $i$ for $2\leq i\leq m$, then fill in
the cells corresponding to $\mathsf{p}(z^{(1)})$ with $1$.
\end{lemma}

\medskip

The collection $\mathcal{P}(m)$ of all the GT\ patterns of $GL_{m}$ with the
function addition forms a semigroup, or more precisely a commutative monoid
with the zero function as its identity, which we shall call the \textit{%
semigroup of patterns} for $GL_{m}$:%
\begin{equation}
\mathcal{P}(m)=\left\{ \mathsf{p}:\Gamma _{m}\rightarrow \mathbb{Z}_{\geq
0}\right\} .  \label{plainpatterns}
\end{equation}%
Then $\mathcal{P}(m)$ is generated by characteristic functions over order
increasing subsets of $\Gamma _{m}$ and these generators correspond to
elements of $L(m,m)$ via Lemma \ref{pattern2tableau}. In this connection,
one can realize the semigroup ring of $\mathcal{P}(m)$ as the Hibi algebra (%
\cite{Hi87}) over the distributive lattice $L(m,m)$. Moreover, the semigroup
structure is compatible with the combinatorial correspondence given in Lemma %
\ref{pattern2tableau} in the following sense: the GT\ pattern corresponding
to a semistandard Young tableau $\mathsf{t}$ or equivalently a multiple
chain $J_{1}\preceq \cdots \preceq J_{u}$ in $L(m,m)$ is the product of the
corresponding GT patterns $\mathsf{p}_{J_{k}}$ in the semigroup $\mathcal{P}%
(m)$, i.e.,%
\begin{equation}
\mathsf{t}=(J_{1}\preceq \cdots \preceq J_{u})\mapsto \mathsf{p}_{\mathsf{t}%
}=\sum\limits_{k=1}^{u}\mathsf{p}_{J_{k}}.  \label{compatible}
\end{equation}%
and also $\mathsf{p}_{J}+\mathsf{p}_{J^{\prime }}=\mathsf{p}_{J\wedge
J^{\prime }}+\mathsf{p}_{J\vee J^{\prime }}$. We refer readers to \cite%
{How05} and \cite{Ki08} for further details.

\subsection{Standard monomials}

Let us review standard monomial theory for the Grassmann variety $Gr(n,m+n)$
of $n$ dimensional spaces in $\mathbb{C}^{m+n}$ and its application to a
presentation of the coordinate ring $\mathbb{C}[M_{n,m}]$ of the space $%
M_{n,m}$. The proofs of the results discussed in this subsection and further
details on the structure of $\mathbb{C}[M_{n,m}]$ can be found in \cite[\S 7]%
{BH93} and \cite[\S 13]{Pr07}.

Note that the Pl\"{u}cker coordinates for $Gr(n,m+n)$ can be matched with
the elementary basis elements of $\bigwedge^{n}\mathbb{C}^{m+n}$. By taking
them as elements of $L(n,m+n)$ with fixed length $n$:%
\begin{equation*}
Pl(n,m+n)=\{K\in L(n,m+n):|K|=n\},
\end{equation*}%
we shall continue to denote by $\delta _{K}\in \mathbb{C}[M_{n,m+n}]$ the
maximum minors over $M_{n,m+n}$ whose columns are indexed by $K\in Pl(n,m+n)$%
. Then for any incomparable pair $K,K^{\prime }\in Pl(n,m+n)$, by applying
the Pl\"{u}cker relations to $\delta _{K}\delta _{K^{\prime }}$, we obtain
its standard expression.

\begin{lemma}
(\cite[p. 234, p. 236]{GL01}) \label{straightening1}i) For $K,K^{\prime }\in
Pl(n,m+n)$, the corresponding product $\delta _{K}\delta _{K^{\prime }}\in 
\mathbb{C}[M_{n,m+n}]$ can be uniquely expressed as a linear combination of
standard monomials, i.e.,%
\begin{equation}
\delta _{K}\delta _{K^{\prime }}=\sum_{r}c_{r}\delta _{T_{r}}\delta
_{T_{r}^{\prime }}  \label{straightening}
\end{equation}%
where $T_{r}\preceq T_{r}^{\prime }$ in $Pl(n,m+n)$ for each $r$.

ii) In the right hand side, $\delta _{K\wedge K^{\prime }}\delta _{K\vee
K^{\prime }}$ appears with coefficient $1$, and $T_{r}\preceq K\wedge
K^{\prime }\preceq K\vee K^{\prime }\preceq T_{r}^{\prime }$ for all $r$.
Moreover, for each $(T_{r},T_{r}^{\prime })\neq (K\wedge K^{\prime },K\vee
K^{\prime })$, let $a$ be the smallest integer such that the sum $s$ of the $%
a$-th entries of $T_{r}$ and $T_{r}^{\prime }$ is different from the sum $%
s_{0}$ of the $a$-th entries of $K$ and $K^{\prime }$. Then $s>s_{0}$.
\end{lemma}

By applying the above identities, one can show that any product $\prod
\delta _{K}\in \mathbb{C}[M_{n,m+n}]$ with $K\in Pl(n,m+n)$ can be uniquely
expressed as a linear combination of standard monomials $\Delta (\mathsf{t}%
)=\prod_{j}\delta _{K_{j}}$ with $K_{1}\preceq K_{2}\preceq \cdots $. See 
\cite{BH93, GL01, Hod43} for further details.

Note that all the maximal minors $\delta _{K}\in \mathbb{C}[M_{n,m+n}]$ are
invariant under the left multiplication of the special linear group $SL_{n}$
on $\mathbb{C}[M_{n,m+n}]$ and in fact the maximal minors generate the
invariant ring $\mathbb{C}[M_{n,m+n}]^{SL_{n}}$. This shows that the
invariant ring $\mathbb{C}[M_{n,m+n}]^{SL_{n}}$ forms an algebra with
straightening laws (ASL) with standard monomials $\Delta (\mathsf{t})$ as
its basis.

\medskip

Now, let us consider the embedding $M_{n,m}\rightarrow M_{n,m+n}$ given by $%
X\longmapsto \left( 
\begin{array}{c|c}
X^{\prime } & I_{n}%
\end{array}%
\right) $ and%
\begin{equation*}
M_{n,m+n}^{0}=\{\left( 
\begin{array}{c|c}
X^{\prime } & I_{n}%
\end{array}%
\right) \}\subset M_{n,m+n}
\end{equation*}%
where $I_{n}$ is the $n\times n$ identity matrix and $X^{\prime }$ is the
matrix obtained by reversing the rows of $X=(x_{i,j})$, i.e., 
\begin{equation*}
X^{\prime }=\left( 
\begin{array}{cccc}
x_{n,1} & \cdots & x_{n,m-1} & x_{n,m} \\ 
x_{n-1,1} &  & x_{n-1,m-1} & x_{n-1,m} \\ 
\vdots & \ddots &  & \vdots \\ 
x_{1,1} & \cdots & x_{1,m-1} & x_{1,m}%
\end{array}%
\right)
\end{equation*}

For each $K\in Pl(n,m+n)$ which is not $[m+1,\cdots ,m+n]$, consider the
restriction $\delta _{K}^{0}$ to $M_{n,m+n}^{0}$ of the maximal minor $%
\delta _{K}\in \mathbb{C}[M_{n,m+n}]$. Then $\delta _{K}^{0}$ is equal to,
up to sign, the minor $\delta _{\xi ^{-1}(K)}\in \mathbb{C}[M_{n,m}]$ over $%
M_{n,m}$ with rows and columns given by the one-line tableau $\xi
^{-1}(K)\in D(n,m)$ via the map $\xi $ defined in Lemma \ref{double2single}.
This provides an algebraic realization of the space $M_{n,m}\cong
M_{n,m+n}^{0}$ as a cell of the Grassmann variety $Gr(n,m+n)$.

\begin{proposition}[{\protect\cite[Lemma 7.2.6]{BH93}\protect\cite[p. 522]%
{Pr07}}]
\label{isom1}The map $\widehat{\xi }$ from $\mathbb{C}[M_{n,m+n}]^{SL_{n}}$
to $\mathbb{C}[M_{n,m}]$:%
\begin{equation*}
\widehat{\xi }:\delta _{K}\mapsto \delta _{\xi ^{-1}(K)}
\end{equation*}%
gives an isomorphism between $\mathbb{C}[M_{n,m+n}]^{SL_{n}}/\ker \widehat{%
\xi }$ and $\mathbb{C}[M_{n,m}]$ where $\ker \widehat{\xi }$ is the ideal
generated by $(\delta _{\lbrack m+1,\cdots ,m+n]}-1)$.
\end{proposition}

Since $[m+1,\cdots ,m+n]$ is the largest element in $Pl(n,m+n)$, standard
monomial basis elements $\prod_{r}\delta _{K_{r}}$ which do not end with $%
\delta _{\lbrack m+1,\cdots ,m+n]}$ form a $\mathbb{C}$-basis for $\mathbb{C}%
[M_{m,n}]^{SL_{n}}/\ker \widehat{\xi }$. Using this map $\widehat{\xi }$, we
can transfer the multiplicative structure of $\mathbb{C}[M_{n,m+n}]^{SL_{n}}$
to $\mathbb{C}[M_{n,m}]$. In particular, since $\xi $ in Lemma \ref%
{double2single} is an order isomorphism, the properties of straightening
laws given in Lemma \ref{straightening1} can be also transferred to any
product $\delta _{\lbrack I:J]}\delta _{\lbrack I^{\prime }:J^{\prime }]}$
in $\mathbb{C}[M_{n,m}]$.

\begin{corollary}
\label{isom2}Let $A$, $B$ be elements in $D(n,m)$. Then in the algebra $%
\mathbb{C}[M_{n,m}]$, the corresponding product $\delta _{A}\delta _{B}$ can
be expressed as a linear combination of standard monomials%
\begin{equation}
\delta _{A}\delta _{B}=\sum_{r}c_{r}\delta _{X_{r}}\delta _{Y_{r}}
\label{straightening22}
\end{equation}%
such that for each $r$ we have either $X_{r}\preceq A\wedge B\preceq A\vee
B\preceq Y_{r}$ or $X_{r}\preceq A\wedge B$ with $\delta _{Y_{r}}=1$, and $%
\delta _{A\wedge B}\delta _{A\vee B}$ appears in the right hand side with
coefficient $1$.
\end{corollary}

Let us consider a non-standard monomial $\prod_{r}\delta _{\lbrack
I_{r}:J_{r}]}$ with $l_{r}=\ell ([I_{r}:J_{r}])$ and $l_{1}\geq l_{2}\geq
\cdots $. Then, by applying the above relations as many times as necessary,
we can express $\prod_{r}\delta _{\lbrack I_{r}:J_{r}]}$ as a linear
combination of standard monomials with shapes $(l_{1}^{\prime
},l_{2}^{\prime },\cdots )^{t} $ such that $\sum_{r=1}^{k}l_{r}^{\prime
}\geq \sum_{r=1}^{k}l_{r}$ for each $k\geq 1$. Hence we can impose a
filtration $\mathsf{F}^{sh}=\{\mathsf{F}_{D}^{sh}\}$ by shapes on the
algebra $\mathbb{C}[M_{n,m}]$, and then consider its associated graded
algebra:%
\begin{eqnarray}
\mathbb{C}[M_{n,m}] &=&\sum\limits_{\ell (D)\leq \min (n,m)}\mathsf{F}%
_{D}^{sh}\left( \mathbb{C}[M_{n,m}]\right)  \label{gr} \\
\mathsf{gr}^{sh}\left( \mathbb{C}[M_{n,m}]\right) &=&\sum\limits_{\ell
(D)\leq \min (n,m)}\mathsf{gr}_{D}^{sh}\left( \mathbb{C}[M_{n,m}]\right) 
\notag
\end{eqnarray}

Finally, we obtain standard monomial theory for $\mathbb{C}[M_{n,m}]$.

\begin{corollary}[{\protect\cite[Theorem 7.2.7]{BH93}\protect\cite[p. 530]%
{Pr07}}]
\label{GLSMT}Standard monomials $\Delta (\mathsf{t})$ form a $\mathbb{C}$%
-basis of $\mathbb{C}[M_{n,m}]$. For its associated graded algebra, the $D$%
-graded component $\mathsf{gr}_{D}^{sh}\left( \mathbb{C}[M_{n,m}]\right) $
is spanned by standard monomials of shape $D$.
\end{corollary}

We refer readers to \cite{BH93, Pr07} and \cite{DEP82, DRS76} for further
details on the ASL structure of $\mathbb{C}[M_{n,m}]$ and double tableaux.

\medskip

\section{Lattice Cone for $M_{n,m}$\label{CMnm}}

A \textit{lattice cone} is the intersection of a convex polyhedral cone in $%
\mathbb{R}^{l}$ for some $l$ with $\mathbb{Z}^{l}$. We shall construct a
lattice cone associated with the space $M_{n,m}$ in terms of order
preserving maps on a subposet of the GT poset $\Gamma _{m+n}$ for $GL_{m+n}$%
. Once $\mathbb{C}[M_{n,m}]$ is shown to be a flat deformation of the Hibi
algebra $\mathcal{H}_{D(n,m)}$\ over $D(n,m)$, the lattice cone for $M_{n,m}$
can be understood as a cone encoding the affine toric variety $Spec(\mathcal{%
H}_{D(n,m)})$. The lattice cone for $M_{n,m}$ is related to the
Gelfand-Tsetlin cones attached to the flag varieties studied in \cite{Ki08,
KM05, MS05}.

\subsection{The semigroup $\mathcal{P}_{n,m}$}

Let us write $\widehat{1}$ for $[m+1,\cdots ,m+n]\in Pl(n,m+n)$. For the
semigroup $\mathcal{P}(m+n)$ defined in (\ref{plainpatterns}), let $\mathcal{%
P}(n,m)$ denote its subsemigroup generated by the GT patterns $\mathsf{p}%
_{J} $ corresponding to $J\in Pl(n,m+n)$. Recall that for a subtractive
subset $B$ of a monoid $(A,+)$, the quotient monoid $A/B$ is defined by the
following equivalence relation: $a\sim a^{\prime }$ if there are $b$ and $%
b^{\prime }$ in $B$ such that $a+b=a^{\prime }+b^{\prime }$.

\begin{definition}
\label{semigroup}The semigroup $\mathcal{P}_{n,m}$ of patterns for $M_{n,m}$
is the quotient of $\mathcal{P}(n,m)$ by the multiples of the pattern $%
\mathsf{p}_{\widehat{1}}$ corresponding to $\widehat{1}$:%
\begin{equation*}
\mathcal{P}_{n,m}=\mathcal{P}(n,m)/\left\langle \mathsf{p}_{\widehat{1}%
}\right\rangle .
\end{equation*}%
Let $\mathbb{C}[\mathcal{P}_{n,m}]$ denote the semigroup ring of $\mathcal{P}%
_{n,m}$.
\end{definition}

We can obtain a more explicit description for the elements of $\mathcal{P}%
_{n,m}$. Simple computations of the bijection given in Lemma \ref%
{pattern2tableau} can easily show the following.

\begin{lemma}
i) For every $\mathsf{p}\in \mathcal{P}(n,m)$, the support of $\mathsf{p}$
is contained in%
\begin{equation*}
\Gamma _{m+n}^{n}=\{z_{j}^{(i)}\in \Gamma _{m+n}:j\leq n\}.
\end{equation*}

ii) If $\mathsf{p}\in \mathcal{P}(n,m)$, then for all $z_{b}^{(a)}\in \Gamma
_{m+n}$ with $z_{b}^{(a)}\geq z_{n}^{(m+n)}$ we have $\mathsf{p}%
(z_{b}^{(a)})=\mathsf{p}(z_{n}^{(m+n)}).$

iii) If $\mathsf{p}=\mathsf{p}_{\widehat{1}}$, then $\mathsf{p}%
(z_{b}^{(a)})=1$\ for $z_{b}^{(a)}\geq z_{n}^{(m+n)}$ and $\mathsf{p}%
(z_{b}^{(a)})=0$\ otherwise.
\end{lemma}

\begin{example}
\label{exGT}The ordered elements $[1,2,3],[1,3,4],[2,5,6],\widehat{1}%
=[5,6,7] $ of $Pl(3,7)$ form the following semistandard tableau $\mathsf{t}$%
\ of shape $(4,4,4)$:%
\begin{equation*}
\young(1125,2356,3467)
\end{equation*}%
and we can visualize its corresponding GT pattern $\mathsf{p}_{\mathsf{t}}$
for $GL_{7}$ by listing its values as%
\begin{equation*}
\begin{array}{ccccccccccccc}
4 &  & 4 &  & 4 &  & 0 &  & 0 &  & 0 &  & 0 \\ 
& 4 &  & 4 &  & 3 &  & 0 &  & 0 &  & 0 &  \\ 
&  & 4 &  & 3 &  & 2 &  & 0 &  & 0 &  &  \\ 
&  &  & 3 &  & 2 &  & 2 &  & 0 &  &  &  \\ 
&  &  &  & 3 &  & 2 &  & 1 &  &  &  &  \\ 
&  &  &  &  & 3 &  & 1 &  &  &  &  &  \\ 
&  &  &  &  &  & 2 &  &  &  &  &  & 
\end{array}%
\end{equation*}%
Note that its support is exactly $\Gamma _{7}^{3}\subset \Gamma _{7}$, and $%
\mathsf{p}(z_{b}^{(a)})=4$ for all $z_{b}^{(a)}\geq z_{3}^{(7)}$.
\end{example}

Then, from the description of supports for $\mathsf{p}\in \mathcal{P}(n,m)$
and $\mathsf{p}_{\widehat{1}}$, we obtain the following characterization for
the elements of the quotient $\mathcal{P}_{n,m}$.

\begin{corollary}
Every element of the semigroup $\mathcal{P}_{n,m}$ of patterns for $M_{n,m}$
can be uniquely represented by an order preserving map to the set of
non-negative integers from 
\begin{equation*}
\Gamma _{n,m}=\Gamma _{m+n}^{n}\backslash \{z_{b}^{(a)}:z_{b}^{(a)}\geq
z_{n}^{(m+n)}\}.
\end{equation*}
\end{corollary}

Then the GT pattern for $GL_{7}$ shown in Example \ref{exGT} corresponds to
the following order preserving map defined on $\Gamma _{3,4}$:%
\begin{equation}
\begin{array}{cccccc}
&  & 3 &  &  &  \\ 
& 3 &  & 2 &  &  \\ 
3 &  & 2 &  & 2 &  \\ 
& 3 &  & 2 &  & 1 \\ 
&  & 3 &  & 1 &  \\ 
&  &  & 2 &  & 
\end{array}
\label{exGT2}
\end{equation}

Now we define the convex polyhedral cone associated with the space $M_{n,m}$
as the collection of all non-negative real valued order preserving maps over 
$\Gamma _{n,m}$:%
\begin{equation*}
\mathcal{C}(M_{n,m})=\{f:\Gamma _{n,m}\rightarrow \mathbb{R}_{\geq 0}\}.
\end{equation*}%
Then, by identifying $f$ with its values $(f(z_{b}^{(a)}))\in \mathbb{R}%
^{nm} $ for $z_{b}^{(a)}\in \Gamma _{n,m}$, we can realize the semigroup $%
\mathcal{P}_{n,m}$ of patterns for $M_{n,m}$ as the set of integral points
in $\mathcal{C}(M_{n,m})$, i.e., the \textit{lattice cone for }$M_{n,m}$:%
\begin{equation*}
\mathcal{P}_{n,m}=\mathcal{C}(M_{n,m})\cap \mathbb{Z}^{nm}.
\end{equation*}

\medskip

\subsection{Degeneration of $\mathbb{C}[M_{n,m}]$}

We want to consider a degeneration of the algebra $\mathbb{C}[M_{n,m}]$. We
shall use basically the same degeneration technique shown in \cite{GL96}.

Recall that the Hibi algebra $\mathcal{H}_{L}$ over a distributive lattice $%
L $ is the quotient ring%
\begin{equation*}
\mathcal{H}_{L}\cong \mathbb{C}[z_{x}:x\in L]/\left( z_{x}z_{y}-z_{x\wedge
y}z_{x\vee y}\right)
\end{equation*}%
of the polynomial ring by the ideal generated by $\left\{
z_{x}z_{y}-z_{x\wedge y}z_{x\vee y}\right\} $. Then the Hibi algebra defines
an affine toric variety (\cite{Hi87}).

To extract the Hibi algebra structure from the algebra $\mathbb{C}[M_{n,m}]$%
, we will use the following weight defined via the correspondence given in
Lemma \ref{double2single}.

\begin{definition}
\label{def weight}Let us fix an integer $N$ greater than $2(n+m)$. For $%
[I:J]\in D(n,m)$ and $\xi ([I:J])=[q_{1},\cdots ,q_{n}]$, we define their
weights $wt([I:J])=wt([q_{1},\cdots ,q_{n}])$ as%
\begin{equation*}
wt([I:J])=\sum_{r\geq 1}(m+r-q_{r})N^{n-r}.
\end{equation*}%
The weight of a double tableau $\mathsf{t}$ consisting of $\left\{
[I_{k}:J_{k}]\right\} $ is defined to be the sum of individual weights,
i.e., $wt(\mathsf{t})=\sum_{k}wt([I_{k}:J_{k}])$.
\end{definition}

We shall assume that, by extending the above formula, the weight of $%
[m+1,\cdots ,m+n]$ is equal to zero. We define the weight of $\Delta (%
\mathsf{t})$ to be the weight of the corresponding double tableau $\mathsf{t}
$.

\begin{proposition}
\label{stmD}For $A,B\in D(n,m)$, let 
\begin{equation*}
\delta _{A}\delta _{B}=\sum_{r}c_{r}\delta _{X_{r}}\delta _{Y_{r}}
\end{equation*}%
be the standard expression of $\delta _{A}\delta _{B}$ given in (\ref%
{straightening22}). Then, $wt(A)+wt(B)\geq wt(X_{r})+wt(Y_{r})$ for all $r$
and the equality holds only for $(X_{r},Y_{r})=(A\wedge B,A\vee B)$.
\end{proposition}

It follows directly from Lemma \ref{straightening1} and Corollary \ref{isom2}%
. Let $\xi (A)=K$, $\xi (B)=K^{\prime }$, $\xi (X_{r})=T_{r}$, and $\xi
(Y_{r})=T_{r}^{\prime }$. If $K=[k_{1},\cdots ,k_{n}]$ and $K^{\prime
}=[k_{1}^{\prime },\cdots ,k_{n}^{\prime }]$, then the $i$-th entry of $%
K\wedge K^{\prime }$ is $\min \{k_{i},k_{i}^{\prime }\}$ and the $i$-the
entry of $K\vee K^{\prime }$ is $\max \{k_{i},k_{i}^{\prime }\}$ for $1\leq
i\leq n$. Therefore, if $(T_{r},T_{r}^{\prime })=(K\wedge K^{\prime },K\vee
K)$, then we have $wt(A)+wt(B)=wt(X_{r})+wt(Y_{r})$. Otherwise, from the
second statement of Lemma \ref{straightening1} we have $%
wt(A)+wt(B)>wt(X_{r})+wt(Y_{r})$. Note that if $\delta _{Y_{r}}$ is $1$,
then as discussed after Definition \ref{def weight} the weight of the
corresponding element $[m+1,\cdots ,m+n]$ is $0$, and therefore we still
have the inequality $wt(A)+wt(B)>wt(X_{r})+wt(Y_{r})=wt(X_{r})$.

\begin{theorem}
The algebra $\mathbb{C}[M_{n,m}]$ is a flat deformation of the Hibi algebra $%
\mathcal{H}_{D(n,m)}$ over $D(n,m)$. More precisely, there is a flat $%
\mathbb{C}[t]$ module whose general fiber is isomorphic to $\mathbb{C}%
[M_{n,m}]$ and special fiber is isomorphic to the Hibi algebra over $D(n,m)$.

\begin{proof}
Let us define a $\mathbb{Z}$-filtration $\mathsf{F}^{wt}=\{\mathsf{F}%
_{d}^{wt}\}$ on $\mathbb{C}[M_{n,m}]$ with respect to the weight $wt$, i.e., 
$\mathsf{F}_{d}^{wt}(\mathbb{C}[M_{n,m}])$ is the $\mathbb{C}$-span of the
set%
\begin{equation*}
\{\Delta (\mathsf{t}):wt(\mathsf{t})\leq d\}.
\end{equation*}%
The filtration $\mathsf{F}^{wt}$ is well defined, since every product $\prod
\delta _{A}$ can be expressed as a linear combination of standard monomials
with smaller weights by the above Proposition. For all pairs $A,B\in D(n,m)$%
, since $wt(A)+wt(B)=wt(A\wedge B)+wt(A\vee B)$, $\delta _{A}\delta _{B}$
and $\delta _{A\wedge B}\delta _{A\vee B}$ belong to the same associated
graded space. Therefore, we have $s_{A}\cdot _{gr}s_{B}=s_{A\wedge B}\cdot
_{gr}s_{A\vee B}$ where $s_{C}$ are elements corresponding to $\delta _{C}$
in the associated graded ring $\mathsf{gr}^{wt}(\mathbb{C}[M_{n,m}])$ of $%
\mathbb{C}[M_{n,m}]$ with respect to the filtration $\mathsf{F}^{wt}$. Then
it is straightforward to show that the associated graded ring $\mathsf{gr}%
^{wt}(\mathbb{C}[M_{n,m}])$ forms the Hibi algebra over $D(n,m)$. From the
general properties of the Rees algebras (e.g., \cite{AB04}), the Rees
algebra $\mathcal{R}^{t}$ of $\mathbb{C}[M_{n,m}]$ with respect to $\mathsf{F%
}^{wt}$:%
\begin{equation*}
\mathcal{R}^{t}=\bigoplus_{d\geq 0}\mathsf{F}_{d}^{wt}(\mathbb{C}%
[M_{n,m}])t^{d}
\end{equation*}%
is flat over $\mathbb{C}[t]$ with its general fiber isomorphic to $\mathbb{C}%
[M_{n,m}]$ and special fiber isomorphic to the associated graded ring which
is $\mathcal{H}_{D(n,m)}$.
\end{proof}
\end{theorem}

\medskip

\subsection{Representations}

Every irreducible polynomial representation of $GL_{k}$ is, via the
correspondence between dominant weights and Young diagrams, uniquely labeled
by a Young diagram with no more than $k$ rows. Let $\rho _{k}^{D}$ denote
the irreducible representation of $GL_{k}$ labeled by Young diagram $D$. Now
we impose an action of $GL_{n}\times GL_{m}$ on the space $M_{n,m}$ by 
\begin{equation}
(g_{1},g_{2})Q=(g_{1}^{t})^{-1}Qg_{2}^{-1}  \label{GLaction}
\end{equation}%
for $g_{1}\in GL_{n}$, $g_{2}\in GL_{m}$, and $Q\in M_{n,m}$. Then we have
the following decomposition of the coordinate ring $\mathbb{C}[M_{n,m}]$
with respect to the action:%
\begin{equation*}
\mathbb{C}[M_{n,m}]=\sum\limits_{\ell (D)\leq \min (n,m)}\rho
_{n}^{D}\otimes \rho _{m}^{D}
\end{equation*}%
where the summation runs over all $D$ of length $\ell (D)$ less than or
equal to $\min (n,m)$. This result is known as $GL_{n}$-$GL_{m}$ duality.
See \cite[Corollary 4.5.19]{GW98} and \cite[Theorem 2.1.2]{How95}.

\medskip

Every minor $\delta _{\lbrack I:J]}$ over $M_{n,m}$ with $[I:J]\in D(n,m)$
is scaled under the action of the diagonal subgroups of $GL_{n}$ and $GL_{m}$%
. Therefore, standard monomials $\Delta (\mathsf{t})=\prod_{k=1}^{r}\delta
_{\lbrack I_{k}:J_{k}]}$ can be seen as joint weight vectors for the
irreducible $GL_{n}\times GL_{m}$ representation $\rho _{n}^{D}\otimes \rho
_{m}^{D}$ where $D$ is equal to $sh(\mathsf{t})$, i.e., $D=(\ell
([I_{1}:J_{1}]),\cdots ,\ell ([I_{r}:J_{r}]))^{t}$. Then, by Corollary \ref%
{GLSMT} the above representation decomposition is compatible with the
associated graded algebra in (\ref{gr}) in that $\mathsf{gr}_{D}^{sh}(%
\mathbb{C}[M_{n,m}])=\rho _{n}^{D}\otimes \rho _{m}^{D}$ and%
\begin{equation*}
\mathsf{gr}^{sh}\left( \mathbb{C}[M_{n,m}]\right) =\sum\limits_{\ell (D)\leq
\min (n,m)}\rho _{n}^{D}\otimes \rho _{m}^{D}.
\end{equation*}

On the other hand, if we write $\mathcal{P}_{n,m}^{D}$ for the collection of
elements $\mathsf{p}$ in the lattice cone $\mathcal{P}_{n,m}$ for $M_{n,m}$\
such that $\mathsf{p}(z^{(m)})=D$, then $\mathcal{P}_{n,m}$ can be expressed
as the following disjoint union:%
\begin{equation*}
\mathcal{P}_{n,m}=\bigcup_{\ell (D)\leq \min (n,m)}\mathcal{P}_{n,m}^{D}
\end{equation*}%
over Young diagrams $D$ of length not more than $\min (n,m)$.

Recall that the GT patterns $\mathsf{p}$ of $GL_{k}$ with fixed top row $%
\mathsf{p}(z^{(k)})=D$ encode weight basis elements for the irreducible
representation $\rho _{k}^{D}$ of $GL_{k}$ with highest weight $D$ (\cite%
{GT50}). Hence the joint weight vectors of $\rho _{n}^{D}\otimes \rho
_{m}^{D}$ can be encoded by pairs of GT patterns of $GL_{n}$ and $GL_{m}$
with the same top row $D$.

\begin{proposition}
\label{cone4GL}The Hibi algebra $\mathcal{H}_{D(n,m)}$ over $D(n,m)$\ is
isomorphic to the semigroup ring $\mathbb{C}[\mathcal{P}_{n,m}]$.
\end{proposition}

\begin{proof}
Since the multiple chains in $D(n,m)$ provide a $\mathbb{C}$-basis of the
Hibi algebra over $D(n,m)$ (\cite{Hi87}), let us find a bijection between
the set of standard tableaux of shape $D$ and $\mathcal{P}_{n,m}^{D}$. For a
standard tableau $\mathsf{t}$ of shape $D$ consisting of $\{[I_{k}:J_{k}]\}$%
, consider the GT pattern $\mathsf{p}$ of $GL_{n+m}$ corresponding to the
semistandard tableau whose columns are $\{\xi ([I_{k}:J_{k}])\}$ where $\xi $
is the bijection given in Lemma \ref{double2single}. Then this
correspondence is injective, and it is straightforward to check that $%
\mathsf{p}$ is an element of $\mathcal{P}(n,m)$ satisfying $\mathsf{p}%
(z^{(m)})=D$. To see surjectivity, note that $\mathsf{p}$ as an element of $%
\mathcal{P}_{n,m}=\mathcal{P}(n,m)/\left\langle \mathsf{p}_{\widehat{1}%
}\right\rangle $\ can be decomposed into two parts having supports in $%
\{z_{b}^{(a)}\in \Gamma _{n,m}:a\geq m\}$ and in $\{z_{b}^{(a)}\in \Gamma
_{n,m}:a\leq m\}$ respectively, and therefore GT patterns of $GL_{n}$ and $%
GL_{m}$ with the same top rows $D$. So we have an one-to-one correspondence
between $\mathcal{P}_{n,m}^{D}$ and the set of standard tableaux of shape $D$%
. This bijection provides an algebra isomorphism from our discussion (\ref%
{compatible}).
\end{proof}

In the proof we used the following observation: $\Gamma _{n,m}$ can be seen
as a gluing of two GT posets $\Gamma _{n}$ and $\Gamma _{m}$ along their top
rows. For instance, the pattern in the quotient $\mathcal{P}_{n,m}$ given in
(\ref{exGT2}) can be seen as a fiber product of GT patterns for $GL_{3}$ and 
$GL_{4}$ over their top rows:%
\begin{equation*}
\begin{array}{ccccc}
3 &  & 2 &  & 2 \\ 
& 3 &  & 2 &  \\ 
&  & 3 &  & 
\end{array}%
\text{ and }%
\begin{array}{ccccccc}
3 &  & 2 &  & 2 &  & 0 \\ 
& 3 &  & 2 &  & 1 &  \\ 
&  & 3 &  & 1 &  &  \\ 
&  &  & 2 &  &  & 
\end{array}%
\end{equation*}

\begin{remark}
\label{gluing}This is a GT pattern version of the correspondence given in
Lemma \ref{double2single} for tableaux. That is, for a standard monomial $%
\Delta (\mathsf{t})=\prod_{k}\delta _{\lbrack I_{k}:J_{k}]}$, let $T^{-}$
and $T^{+}$ be the semistandard tableaux whose columns are $\{I_{k}\}$ and $%
\{J_{k}\}$ respectively, and let $\xi ([I_{k}:J_{k}])=K_{k}\in Pl(n,m+n)$
for each $k$. Then the pattern $\mathsf{p}\in \mathcal{P}(n+m)$
corresponding to the semistandard tableau with columns $\{K_{k}\}$ can be,
as an element of $\mathcal{P}_{n,m}$, identified with the gluing of $\mathsf{%
p}_{-}$ and $\mathsf{p}_{+}$ along their top rows where $\mathsf{p}_{-}\in 
\mathcal{P}(n)$ and $\mathsf{p}_{+}\in \mathcal{P}(m)$ are the GT\ patterns
corresponding to $T^{-}$ and $T^{+}$ respectively.
\end{remark}

\medskip

\section{Standard Monomial Theory for $\mathcal{N}_{k,2n}$\label{Nk2n}}

In this section, we define the nullcone $\mathcal{N}_{k,2n}$ in the
multi-vector representation of the symplectic group and consider the $%
GL_{k}\times Sp_{2n}$ action on it. Then we study standard monomial theory
and a toric degeneration of $\mathcal{N}_{k,2n}$. Having an explicit
description of the standard monomials for $\mathbb{C}[M_{k,2n}]$, we develop
a relative theory to $M_{k,2n}$ for $\mathcal{N}_{k,2n}$ by investigating
the defining ideal of $\mathcal{N}_{k,2n} $.

\subsection{Nullcone for $Sp_{2n}$}

For the space $\mathbb{C}^{2n}$\ with the elementary basis $\{e_{i}\}$, let
us fix our skew symmetric bilinear form $\left\langle ,\right\rangle $ on it
such that for every $i$, $e_{2i-1}$ and $e_{2i}$ form an isotropic pair with 
$\left\langle e_{2i-1},e_{2i}\right\rangle =1$. We can consider the space $%
M_{k,2n}$ of $k\times 2n$ complex matrices as $k$ copies of $\mathbb{C}^{2n}$
with the natural action of the symplectic group $Sp_{2n}$. Then by the first
fundamental theorem of invariant theory (e.g., \cite[Theorem 4.2.2]{GW98}%
\cite[Theorem 3.8.3.2]{How95}), the $Sp_{2n}$ invariants of $\mathbb{C}%
[M_{k,2n}]$ are generated by the basic invariants $r_{ij}=\left\langle
v_{i},v_{j}\right\rangle $ obtained from row vectors $v_{i}$ and $v_{j}$, or
in terms of the coordinates specified in (\ref{coordinates}),%
\begin{equation*}
r_{ij}=\sum_{u=1}^{n}\left( x_{i,2u-1}x_{j,2u}-x_{j,2u-1}x_{i,2u}\right)
\end{equation*}%
for $1\leq i<j\leq k$.

\begin{definition}
The\textit{\ nullcone }$\mathcal{N}_{k,2n}$\textit{\ for} $Sp_{2n}$ is the
subvariety of $M_{k,2n}$\ defined by the $Sp_{2n}$-invariants with vanishing
constant terms.
\end{definition}

If we let $\mathcal{I}$ denote the ideal of $\mathbb{C}[M_{k,2n}]$ generated
by $\{r_{ij}:1\leq i<j\leq k\}$, then it is a radical ideal and the
coordinate ring of $\mathcal{N}_{k,2n}$ is%
\begin{equation*}
\mathcal{R}(\mathcal{N}_{k,2n})=\mathbb{C}[M_{k,2n}]/\mathcal{I}
\end{equation*}%
See \cite[Theorem 3.8.6.2]{How95}. One can also study the nullcone $\mathcal{%
N}_{k,2n}$ as the zero fiber $\pi ^{-1}(0)$ of the quotient $\pi
:M_{k,2n}\rightarrow M_{k,2n}//Sp_{2n}$ and investigate the orbit structure.
See \cite{Kr85} for this direction.

\medskip

From the action of $GL_{k}\times GL_{2n}$ on $M_{k,2n}$ given in (\ref%
{GLaction}), by taking $Sp_{2n}$ as a subgroup of $GL_{2n}$, we can consider
the action of $GL_{k}\times Sp_{2n}$ on $\mathbb{C}[M_{k,2n}]$. Moreover,
since $GL_{k}$ and $Sp_{2n}$ commute with each other in this action, the
ideal $\mathcal{I}$ is stable under $GL_{k}\times Sp_{2n}$. Therefore, we
can regard $\mathcal{R}(\mathcal{N}_{k,2n})$ as a $GL_{k}\times Sp_{2n}$
stable complement of $\mathcal{I}$.

Recall that by highest weight theory, every polynomial representation of $%
GL_{k}$ and $Sp_{2n}$ can be uniquely labeled by a Young diagram with no
more than $k$ and $n$ rows respectively. We let $\rho _{k}^{D}$ and $\sigma
_{2n}^{D}$ denote the irreducible representations of $GL_{k}$ and $Sp_{2n}$
respectively labeled by Young diagram $D$.

\begin{proposition}[{\protect\cite[Theorem 3.8.6.2]{How95}}]
Under the action of $GL_{k}\times Sp_{2n}$, we have the following
decomposition:%
\begin{equation*}
\mathcal{R}(\mathcal{N}_{k,2n})=\sum_{r(D)\leq \min (n,k)}\rho
_{k}^{D}\otimes \sigma _{2n}^{D}
\end{equation*}%
where the summation runs over all Young diagrams $D$ with length no more
than $\min (k,n)$.
\end{proposition}

We remark that the space $\mathcal{H}(M_{k,2n})$ of $Sp_{2n}$-\textit{%
harmonics} in $\mathbb{C}[M_{k,2n}]$ can be defined by the kernel of the
symplectic analogs of Laplacian differential operators. Then, as is the case
for $\mathcal{R}(\mathcal{N}_{k,2n})$, the space of harmonics is stable
under the action of $GL_{k}\times Sp_{2n}$. In fact, $\mathcal{H}(M_{k,2n})$
and $\mathcal{R}(\mathcal{N}_{k,2n})$ share the same decomposition under the
action of $GL_{k}\times Sp_{2n}$ (\cite{GW98, How95}). Therefore, our
results may be used to study the space of harmonics.

\subsection{Standard Monomials for $\mathcal{N}_{k,2n}$}

Let us fix some notations. We write $\omega $ for the following $Sp_{2n}$%
-invariant element in $\bigwedge\nolimits^{2}\mathbb{C}^{2n}$:%
\begin{equation*}
\omega =\sum\nolimits_{u=1}^{n}e_{2u-1}\wedge e_{2u}\in
\bigwedge\nolimits^{2}\mathbb{C}^{2n}.
\end{equation*}%
For $J=[j_{1},\cdots ,j_{p}]\in L(2n,2n)$, write $e_{J}$ for the elementary
basis element $e_{j_{1}}\wedge e_{j_{2}}\wedge \cdots \wedge e_{j_{p}}\in
\bigwedge\nolimits^{p}\mathbb{C}^{2n}$.

\begin{definition}
An $\omega $-sum of $Sp_{2n}$ is a linear combination $\sum%
\nolimits_{d=1}^{r}c_{d}J_{d}$ of elements from $L(2n,2n)$ such that%
\begin{equation*}
\sum\nolimits_{d=1}^{r}c_{d}e_{J_{d}}\in \omega \wedge
(\bigwedge\nolimits^{p-2}\mathbb{C}^{2n})
\end{equation*}%
for some $p\geq 2$. We denote by $\Omega _{2n}$ the collection of $\omega $%
-sums of $Sp_{2n}$.
\end{definition}

\begin{proposition}
The following set generates the ideal $\mathcal{I}\subset \mathbb{C}%
[M_{k,2n}]$ of the nullcone $\mathcal{N}_{k,2n}$:%
\begin{equation*}
\Theta =\{\sum\nolimits_{d}c_{d}\delta _{\lbrack
I:J_{d}]}:\sum\nolimits_{d}c_{d}J_{d}\in \Omega _{2n}\}.
\end{equation*}
\end{proposition}

\begin{proof}
The ideal generated by$\mathcal{\ }\Theta $ contains $\mathcal{I}$, because
the basic $Sp_{2n}$-invariants $r_{ij}$ are elements of $\Theta $ with $%
I=[i,j]$, $J_{d}=[2d-1,2d]$ and $c_{d}=1$ for $1\leq d\leq n$. For $%
e_{k_{1}}\wedge e_{k_{2}}\wedge \cdots \wedge e_{k_{p-2}}\in
\bigwedge\nolimits^{p-2}\mathbb{C}^{2n}$, let us consider the following
elements in $\omega \wedge (\bigwedge\nolimits^{p-2}\mathbb{C}^{2n})$:%
\begin{eqnarray*}
\omega \wedge \left( e_{k_{1}}\wedge e_{k_{2}}\wedge \cdots \wedge
e_{k_{p-2}}\right) &=&\sum\nolimits_{u=1}^{n}e_{2u-1}\wedge e_{2u}\wedge
e_{k_{1}}\wedge e_{k_{2}}\wedge \cdots \wedge e_{k_{p-2}} \\
&=&\sum\nolimits_{u=1}^{n}\sigma _{u}\left( e_{j_{1}}\wedge e_{j_{2}}\wedge
\cdots \wedge e_{j_{p}}\right) \\
&=&\sum\nolimits_{u=1}^{n}\sigma _{u}e_{J_{u}}
\end{eqnarray*}%
where $\{j_{1},\cdots ,j_{p}\}=\{2u-1,2u,k_{1},\cdots ,k_{p-2}\}$ with $%
j_{1}\leq \cdots \leq j_{p}$. If there is a repetition in $%
\{2u-1,2u,k_{1},\cdots ,k_{p-2}\}$, then $\sigma _{u}=0$. If there is no
repetition in $\{2u-1,2u,k_{1},\cdots ,k_{p-2}\}$, then $\sigma _{u}$ is the
signature of the permutation sorting $2u-1,2u,k_{1},\cdots ,k_{p-2}$ in
increasing order. Since $\omega \wedge (\bigwedge^{p-2}\mathbb{C}^{2n})$ is
spanned by these elements, the elements of $\Theta $ are linear combinations
of their associated elements $\sum\nolimits_{u=1}^{n}\sigma _{u}\delta
_{\lbrack I:J_{u}]}$. The column expansions for the determinants $\delta
_{\lbrack I:J_{u}]}$ show that $\sum\nolimits_{u}\sigma _{u}\delta _{\lbrack
I:J_{u}]}$ is an element of the ideal generated by the basic $Sp_{2n}$%
-invariants $\{r_{ij}\}$. This shows that $\Theta $ is contained in$\mathcal{%
\ I}$, and therefore $\Theta $ generates the ideal$\mathcal{\ I}$.
\end{proof}

Next, we characterize standard monomials of $\mathbb{C}[M_{k,2n}]$
associated with the elements of the ideal $\mathcal{I}$, and then we define
standard monomials for the quotient $\mathcal{R}(\mathcal{N}_{k,2n})$. Let
us impose the lexicographic order on the elements of the same length in $%
L(2n,2n)$. We say $[i_{1},\cdots ,i_{p}]>_{lex}[j_{1},\cdots ,j_{p}]$ if the
left-most nonzero entry of $(i_{1}-j_{1},\cdots ,i_{p}-j_{p})\in \mathbb{Z}%
^{p}$ is positive. We fix the element $\widetilde{J}=[1,3,\cdots ,2n-1]$ of
length $n$ having $2d-1$ as its $d$-th smallest entry for $1\leq d\leq n$.

\begin{lemma}
\label{Spstratightening}Let $\sum\nolimits_{d=1}c_{d}J_{d}$ be an $\omega $%
-sum of $Sp_{2n}$. Then the smallest element $J_{1}$ among $\{J_{d}\}$ with
respect to the lexicographic order satisfies $J_{1}\nsucceq \widetilde{J}$.
Conversely, if $J_{1}\in L(2n,2n)$ satisfies $J_{1}\nsucceq \widetilde{J}$,
there is an $\omega $-sum of $Sp_{2n}$ whose smallest non-zero term with
respect to the lexicographic order is $J_{1}$.
\end{lemma}

In particular, note that if $\ell (J)>n$, then $J\nsucceq \widetilde{J}$.
The above Lemma is from computations of the fundamental representations of $%
Sp_{2n}$, which can be realized in the quotient of $\bigwedge \mathbb{C}%
^{2n} $ by the ideal generated by $\omega $. Its proof can be found in \cite[%
Proposition 5.6, Proposition 5.9]{Ki08} and \cite[\S 17]{FH91}. See also 
\cite{Be86} for a combinatorial description of such computations.

\begin{definition}
Let us define a distributive lattice $D(\mathcal{N})$ for $\mathcal{N}=%
\mathcal{N}_{k,2n}$ as%
\begin{equation*}
D(\mathcal{N})=\{[I:J]\in D(k,2n):\ell ([I:J])\leq \min (k,n)\text{ and }%
J\succeq \widetilde{J}\}.
\end{equation*}%
A multiple chain $\mathsf{t}=(X_{1}\preceq X_{2}\preceq \cdots )$ in the
poset $D(\mathcal{N})$ is called an $\mathcal{N}$-standard tableau, and the
corresponding monomial $\Delta (\mathsf{t})=\prod_{r}\delta _{X_{r}}\in 
\mathbb{C}[M_{k,2n}]$ is called an $\mathcal{N}$-standard monomial.
\end{definition}

\begin{proposition}
\label{Nwt}i) For $A,B\in D(\mathcal{N})$, the corresponding product $\delta
_{A}\delta _{B}$ in $\mathcal{R}(\mathcal{N}_{k,2n})$ can be expressed as a
linear combination of $\mathcal{N}$-standard monomials:%
\begin{equation*}
\delta _{A}\delta _{B}=\sum_{r}c_{r}\delta _{X_{r}}\delta _{Y_{r}}
\end{equation*}%
where $\delta _{Y_{r}}$ can possibly be $1$, and $\delta _{A\wedge B}\delta
_{A\vee B}$ appears in the right hand side with coefficient $1$.

ii) Moreover, in the above expression, $wt(A)+wt(B)\geq wt(X_{r})+wt(Y_{r})$
and the equality holds only for $(X_{r},Y_{r})=(A\wedge B,A\vee B)$.

\begin{proof}
Note that for $A,B\in D(\mathcal{N})$, $A\wedge B$ and $A\vee B$ belong to $%
D(\mathcal{N})$ and the corresponding element $\delta _{A\wedge B}\delta
_{A\vee B}$ appears in the standard expression of $\delta _{A}\delta _{B}$
in $\mathbb{C}[M_{k,2n}]$ by Proposition \ref{stmD}. Then the first
statement follows easily from the following computation: starting from the
standard expression $\sum_{r}c_{r}\delta _{X_{r}}\delta _{Y_{r}}$ of $\delta
_{A}\delta _{B}$ in $\mathbb{C}[M_{k,2n}]$ given in Proposition \ref{stmD},
if there is $X_{r}=[I:J]$ which is in $D(k,2n)\backslash D(\mathcal{N})$
then we can obtain the $\mathcal{N}$-standard expression of $\delta
_{X_{r}}\delta _{Y_{r}}$ by successive applications of the elements $\left(
\delta _{\lbrack I:J]}-\sum_{d}s_{d}\delta _{\lbrack I:J_{d}]}\right) $ in
the generating set $\Theta $ of the ideal $\mathcal{I}$ such that $%
[I:J_{d}]>_{lex}[I:J]$ for all $d$, combined with the relations in
Proposition \ref{stmD} if necessary. We can always find such elements in $%
\Theta $\ by Lemma \ref{Spstratightening}. For the second statement, note
that after replacing a non $\mathcal{N}$-standard term $\delta
_{X_{r}}\delta _{Y_{r}}$ by $(\sum_{d}s_{d}\delta _{\lbrack I:J_{d}]})\delta
_{Y_{r}}$, the weights of new terms $wt([I:J_{d}])+wt(Y_{r})$ are strictly
smaller than $wt(X_{r})+wt(Y_{r})$. If a term $\delta _{X_{r}}\delta
_{Y_{r}} $ is already $\mathcal{N}$-standard, then the inequality of the
weight directly follows from Proposition \ref{stmD}.
\end{proof}
\end{proposition}

\medskip

Now we state standard monomial theory for the nullcone $\mathcal{N}_{k,2n}$
and show its degeneration by the same method used for $\mathbb{C}[M_{n,m}]$
in Section \ref{CMnm}. 

\begin{theorem}
\label{standardmonomialN}The $\mathcal{N}$-standard monomials$\ $in $\mathbb{%
C}[M_{k,2n}]$ project to a $\mathbb{C}$-basis of $\mathcal{R}(\mathcal{N}%
_{k,2n})$. In particular, the $\mathcal{N}$-standard monomials $\Delta (%
\mathsf{t})$ of shape $sh(\mathsf{t})=D$ project to a basis of $\rho
_{k}^{D}\otimes \sigma _{2n}^{D}$.

\begin{proof}
For a standard monomial $\prod_{r}\delta _{\lbrack I_{r}:J_{r}]}$\ of $%
\mathbb{C}[M_{k,2n}]$, if $J_{s}\nsucceq \widetilde{J}$ for some $s$, then
we find $\delta _{\lbrack I_{s}:J_{s}]}-\sum_{d}c_{d}\delta _{\lbrack
I_{s}:J_{d,s}]}\in \Theta $ with $J_{d,s}>_{lex}J_{s}$ for all $d$ by Lemma %
\ref{Spstratightening}. Then 
\begin{equation*}
\left( \delta _{\lbrack I_{s}:J_{s}]}-\sum_{d}c_{d}\delta _{\lbrack
I_{s}:J_{d,s}]}\right) \prod_{r\neq s}\delta _{\lbrack I_{r}:J_{r}]}
\end{equation*}%
is in the ideal $\mathcal{I}$ having the monomial $\prod_{r}\delta _{\lbrack
I_{r}:J_{r}]}$ as its initial term. By repeating this procedure, combined
with the relation in Proposition \ref{Nwt}\ if necessary, we can express $%
\prod_{r}\delta _{\lbrack I_{r}:J_{r}]}$ as a linear combination of $%
\mathcal{N}$-standard monomials. This implies that $\mathcal{N}$-standard
monomials project to a spanning set of the space $\mathcal{R}(\mathcal{N}%
_{k,2n})$. Moreover, since the subspace of $\mathcal{R}(\mathcal{N}_{k,2n})$
spanned by $\mathcal{N}$-standard monomials of shape $D$ is stable under $%
GL_{k}\times Sp_{2n}$, the dimension of the space $\rho _{k}^{D}\otimes
\sigma _{2n}^{D}$ is less than or equal to the number of $\mathcal{N}$%
-standard monomials of shape $D$.

Now we claim that the number of $\mathcal{N}$-standard monomials of shape $D$
is exactly the dimension of the space $\rho _{k}^{D}\otimes \sigma _{2n}^{D}$%
. For an $\mathcal{N}$-standard monomial $\prod_{r}\delta _{\lbrack
I_{r}:J_{r}]}$ of shape $D$, the row indices $\{I_{r}\}$ form a semistandard
tableau $T^{-}$ of shape $D$ with entries from $\{1,\cdots ,k\}$. Then, the
number of such semistandard tableaux is equal to the dimension of $\rho
_{k}^{D}$ (e.g., \cite{FH91, GW98}). On the other hand, the column indices $%
\{J_{r}\}$ form a semistandard tableau $T^{+}$ of shape $D$ such that each
column is greater than or equal to $\widetilde{J}$. The set of all possible
such semistandard tableaux $T^{+}$ with entries from $\{1,\cdots ,2n\}$
labels the weight basis of $\sigma _{2n}^{D}$ (e.g., \cite{Be86, Ki08}).
Therefore, the number of all the $\mathcal{N}$-standard monomials $\Delta (%
\mathsf{t})$ with $sh(\mathsf{t})=D$\ is equal to the dimension of $\rho
_{k}^{D}\otimes \sigma _{2n}^{D}$. This finally shows that the $\mathcal{N}$%
-standard monomials with shape $D$ project to a $\mathbb{C}$-basis of $\rho
_{k}^{D}\otimes \sigma _{2n}^{D}$.
\end{proof}
\end{theorem}

\begin{theorem}
The algebra $\mathcal{R}(\mathcal{N}_{k,2n})$ is a flat deformation of the
Hibi algebra over $D(\mathcal{N})$. More precisely, there is a flat $\mathbb{%
C}[t]$ module whose general fiber is $\mathcal{R}(\mathcal{N}_{k,2n})$ and
special fiber is isomorphic to the Hibi algebra $\mathcal{H}_{D(\mathcal{N}%
)} $ over $D(\mathcal{N})$.
\end{theorem}

\begin{proof}
From Theorem \ref{standardmonomialN}, every element of $\mathcal{R}(\mathcal{%
N}_{k,2n})$ can be uniquely expressed as a linear combination of $\mathcal{N}
$-standard monomials $\Delta (\mathsf{t})$. Hence, we can impose the same
filtration $\mathsf{F}^{wt}$ of $\mathbb{C}[M_{k,2n}]$ on $\mathcal{R}(%
\mathcal{N}_{k,2n})$ via the weights $wt$ on $\mathcal{N}$-standard
monomials (Definition \ref{def weight} with $D(k,2n)$): $\mathsf{F}_{d}^{wt}(%
\mathcal{R}(\mathcal{N}_{k,2n}))$ is the $\mathbb{C}$-span of the set 
\begin{equation*}
\{\Delta (\mathsf{t}):wt(\mathsf{t})\leq d\}.
\end{equation*}%
This filtration is well defined, since in the standard expression $%
\sum_{r}c_{r}\Delta (\mathsf{t}_{r})$ of any product $\prod \delta _{A}$,
the weights $wt(\mathsf{t}_{r})$ are smaller than the weight of $\prod
\delta _{A}$ by Proposition \ref{Nwt}. Moreover, since the equality holds
only for $(X_{r},Y_{r})=(A\wedge B,A\vee B)$, as in the case for the space $%
M_{k,2n}$, we have the relation $s_{A}\cdot _{gr}s_{B}=s_{A\wedge B}\cdot
_{gr}s_{A\vee B}$ where $s_{C}$ are elements corresponding to $\delta _{C}$
in the associated graded algebra $\mathsf{gr}^{wt}(\mathcal{R}(\mathcal{N}%
_{k,2n}))$ with respect to $\mathsf{F}^{wt}$. Therefore it is easy to see
that the associated graded algebra forms the Hibi algebra over $D(\mathcal{N}%
)$. Now for the flat degeneration, we can construct the Rees algebra $%
\mathcal{R}^{t}$:%
\begin{equation*}
\mathcal{R}^{t}=\bigoplus_{d\geq 0}\mathsf{F}_{d}^{wt}(\mathcal{R}(\mathcal{N%
}_{k,2n}))t^{d}
\end{equation*}%
with respect to $\mathsf{F}^{wt}$, then from the general properties of the
Rees algebras (e.g., \cite{AB04}), $\mathcal{R}^{t}$ is flat over $\mathbb{C}%
[t]$ with general fiber isomorphic to $\mathcal{R}(\mathcal{N}_{k,2n})$ and
special fiber isomorphic to the associated graded algebra which is $\mathcal{%
H}_{D(\mathcal{N})}$.
\end{proof}

\medskip

\section{Lattice Cone for $\mathcal{N}_{k,2n}$\label{CNk2n}}

In this section we study a lattice cone associated with $\mathcal{N}_{k,2n}$%
. As is the case for $M_{n,m}$, it turns out that each point in the lattice
cone for $\mathcal{N}_{k,2n}$ can be identified with a pair of
Gelfand-Tsetlin patterns.

\begin{proposition}
\label{support}Let $\Delta (\mathsf{t})=\prod_{r}\delta _{\lbrack
I_{r}:J_{r}]}$ be an $\mathcal{N}$-standard monomial and $\mathsf{p}_{%
\mathsf{t}}\in \mathcal{P}_{k,2n}$ be the pattern corresponding to $\Delta (%
\mathsf{t})$. Then, $\mathsf{p}_{\mathsf{t}}$ has its support in the
following subposet of $\Gamma _{k,2n}^{{}}$:%
\begin{equation}
\digamma _{k,2n}=\Gamma _{k,2n}^{{}}\backslash \left( A\cup B\right)
\label{N poset}
\end{equation}%
where the subsets $A$ and $B$ of $\Gamma _{2n+k}$ are defined as%
\begin{eqnarray*}
A &=&\{z_{b}^{(a)}\in \Gamma _{2n+k}:a\leq 2n\text{ and }b>\left( a+1\right)
/2\}; \\
B &=&\{z_{b}^{(a)}\in \Gamma _{2n+k}:z_{\min (k,n)+1}^{(2n)}\geq
z_{b}^{(a)}\}.
\end{eqnarray*}
\end{proposition}

We shall prove it in a few steps. First, note that if a GT pattern $\mathsf{p%
}\in \mathcal{P}(k+2n)$ of $GL_{2n+k}$ corresponds to an $\mathcal{N}$%
-standard monomial, then this Proposition says that the length of the $2n$%
-th row is at most $\min (k,n)$ and the support of $\mathsf{p}$
corresponding to the bottom $2n$ rows is contained in the \textquotedblleft
left half\textquotedblright\ of $\Gamma _{2n}\subset \Gamma _{2n+k}$.

\begin{example}
\label{exampleSp}For $k=4$ and $n=3$, let us consider the following $%
\mathcal{N}$-standard monomial for $\mathcal{N}_{4,6}$:%
\begin{equation*}
\delta _{\lbrack 123:135]}\delta _{\lbrack 124:136]}\delta _{\lbrack
12:24]}\delta _{\lbrack 13:35]}\delta _{\lbrack 1:4]}
\end{equation*}%
Then the corresponding semistandard tableau with respect to $\xi $ given in
Lemma \ref{double2single} is the following chain in $Pl(4,10)$:%
\begin{equation*}
\young(11234,33457,56778,78899)
\end{equation*}%
and we can visualize its corresponding GT pattern for $GL_{10}$ by listing
its function values as%
\begin{equation*}
\begin{array}{ccccccccccccccccccc}
5 &  & 5 &  & 5 &  & 5 &  & 0 &  & 0 &  & 0 &  & 0 &  & 0 &  & 0 \\ 
& 5 &  & 5 &  & 5 &  & 5 &  & 0 &  & 0 &  & 0 &  & 0 &  & 0 &  \\ 
&  & 5 &  & 5 &  & 5 &  & 3 &  & 0 &  & 0 &  & 0 &  & 0 &  &  \\ 
&  &  & 5 &  & 5 &  & 4 &  & 1 &  & 0 &  & 0 &  & 0 &  &  &  \\ 
&  &  &  & 5 &  & 4 &  & 2 &  & 0 &  & 0 &  & 0 &  &  &  &  \\ 
&  &  &  &  & 5 &  & 4 &  & 1 &  & 0 &  & 0 &  &  &  &  &  \\ 
&  &  &  &  &  & 5 &  & 3 &  & 0 &  & 0 &  &  &  &  &  &  \\ 
&  &  &  &  &  &  & 4 &  & 2 &  & 0 &  &  &  &  &  &  &  \\ 
&  &  &  &  &  &  &  & 3 &  & 0 &  &  &  &  &  &  &  &  \\ 
&  &  &  &  &  &  &  &  & 2 &  &  &  &  &  &  &  &  & 
\end{array}%
\end{equation*}%
Note that the length of the $6$-th row $(5,4,2)$ is $3$ and the support
corresponding to the bottom six rows is contained only in the
\textquotedblleft left half\textquotedblright\ of the poset of $\Gamma
_{6}\subset \Gamma _{10}$.
\end{example}

Let us write $\widehat{1}$ for $[2n+1,\cdots ,2n+k]\in Pl(k,2n+k)$, and
write $\mathsf{p}_{\widehat{1}}$ for the GT pattern of $GL_{2n+k}$
corresponding to $\widehat{1}$ via Lemma \ref{pattern2tableau}. Recall that $%
\widetilde{J}=[1,3,\cdots ,2n-1]\in L(2n,2n)$.

\begin{lemma}
Let $[I:J]\in D(k,2n)$ be an one-line tableau and $K=\xi ([I:J])$ be the
element of $Pl(k,2n+k)$ corresponding to $[I:J]$ via Lemma \ref%
{double2single}. Let $supp(\mathsf{p}_{K})$ denote the support of the GT
pattern $\mathsf{p}_{K}$\ corresponding to $K$. If $J\nsucceq \widetilde{J}$%
, then the following intersection%
\begin{equation*}
supp(\mathsf{p}_{K})\cap \{z_{b}^{(a)}\in \Gamma _{2n+k}:a\leq 2n\text{ and }%
b>\left( a+1\right) /2\}
\end{equation*}%
is non-empty. Conversely, if $J\succeq \widetilde{J}$, then the intersection
is empty.
\end{lemma}

\begin{proof}
This is an easy computation similar to \cite[Lemma 5.11]{Ki08}.
\end{proof}


Recall that for each standard monomial $\Delta (\mathsf{t})=\prod_{r}\delta
_{\lbrack I_{r}:J_{r}]}$ of $\mathbb{C}[M_{k,2n}]$, by the bijection
constructed in the proof of Proposition \ref{cone4GL}, we can find its
corresponding pattern $\mathsf{p}_{\mathsf{t}}\in \mathcal{P}_{k,2n}$ for $%
M_{k,2n}$. More precisely, $\mathsf{p}_{\mathsf{t}}$ as an element of $%
\mathcal{P}(k,2n)$ is the sum of GT patterns $\mathsf{p}_{r}$ of $GL_{2n+k}$
corresponding to $\xi ([I_{r}:J_{r}])$ where $\xi $ is the bijection given
in Lemma \ref{double2single}. If $\Delta (\mathsf{t})=\prod_{r}\delta
_{\lbrack I_{r}:J_{r}]}$ is a $\mathcal{N}$-standard monomial, then $%
J_{r}\succeq \widetilde{J}$ for all $r$. Therefore, by the above Lemma, the
support of $\mathsf{p}_{\mathsf{t}}\in \mathcal{P}_{k,2n}$ does not
intersect with $\{z_{b}^{(a)}\in \Gamma _{k+2n}:a\leq 2n$ and $b>\left(
a+1\right) /2\}$. Also, note that $J\nsucceq \widetilde{J}$ if $\ell (J)>n$
and that $\ell (I)\leq k$. Then the condition that $\mathsf{p}_{\mathsf{t}}$
is supported in $\Gamma _{k,2n}\backslash \{z_{b}^{(a)}\in \Gamma
_{k+2n}:z_{\min (k,n)+1}^{(2n)}\geq z_{b}^{(a)}\}$ follows from the fact $%
\ell ([I:J])\leq \min (k,n)$ for $[I:J]\in D(\mathcal{N})$. This finishes
the proof of Proposition \ref{support}.

\medskip

Now by using the poset identified in (\ref{N poset}), we can define the
semigroup and the cone for $\mathcal{N}_{k,2n}$.

\begin{definition}
The semigroup $\mathcal{P}(\mathcal{N}_{k,2n})$ of patterns for $\mathcal{N}%
_{k,2n}$ is the set of order preserving maps from $\digamma _{k,2n}$ to the
set of non-negative integers with the usual function addition as its
product. We let $\mathbb{C}[\mathcal{P}(\mathcal{N}_{k,2n})]$ denote the
semigroup ring of $\mathcal{P}(\mathcal{N}_{k,2n})$.
\end{definition}

We can define the convex polyhedral cone associated with the space $\mathcal{%
N}_{k,2n}$ as the collection of all non-negative real valued order
preserving maps on $\digamma _{k,2n}$:%
\begin{equation*}
\mathcal{C}(\mathcal{N}_{k,2n})=\{f:\digamma _{k,2n}\rightarrow \mathbb{R}%
_{\geq 0}\}.
\end{equation*}%
Then, by identifying $f$ with its values $(f(z_{b}^{(a)}))\in \mathbb{R}^{N}$
for $z_{b}^{(a)}\in \digamma _{k,2n}$, we can realize the semigroup $%
\mathcal{P}(\mathcal{N}_{k,2n})$ of patterns for $\mathcal{N}_{k,2n}$ as the
intersection of $\mathcal{C}(\mathcal{N}_{k,2n})$ with $\mathbb{Z}^{N}$,
i.e., the \textit{lattice cone for }$N_{k,2n}$:%
\begin{equation*}
\mathcal{P}(\mathcal{N}_{k,2n})=\mathcal{C}(\mathcal{N}_{k,2n})\cap \mathbb{Z%
}^{N}
\end{equation*}%
where $N$ is equal to the number of elements in the poset $\digamma _{k,2n}$%
. Let us denote by $\mathcal{P}(\mathcal{N}_{k,2n})_{D}$ the collection of $%
\mathsf{p}\in \mathcal{P}(\mathcal{N}_{k,2n})$ whose $2n$-th row is equal to
Young diagram $D$, i.e., $\mathsf{p}(z^{(2n)})=D$. Then the lattice cone for 
$\mathcal{N}_{k,2n}$ can be expressed as the disjoint union%
\begin{equation*}
\mathcal{P}(\mathcal{N}_{k,2n})=\bigcup_{D}\mathcal{P}(\mathcal{N}%
_{k,2n})_{D}
\end{equation*}%
over all $D$ with $\ell (D)\leq \min (k,n)$.

The following is a lattice cone version of Theorem \ref{standardmonomialN}.

\begin{proposition}
\label{cone4Sp}i) For $\mathcal{N}=\mathcal{N}_{k,2n}$, the Hibi algebra $%
\mathcal{H}_{D(\mathcal{N})}$ over $D(\mathcal{N})$ is isomorphic to the
semigroup ring $\mathbb{C}[\mathcal{P}(\mathcal{N}_{k,2n})]$. ii) There is
an one-to-one correspondence between $\mathcal{P}(\mathcal{N}_{k,2n})_{D}$
and the set of weight vectors for $\rho _{k}^{D}\otimes \sigma _{2n}^{D}$ in 
$\mathcal{R}(\mathcal{N}_{k,2n})$.
\end{proposition}

\begin{proof}
The proof of the first statement is parallel to that of Proposition \ref%
{cone4GL}. Note that the multiple chains in $D(\mathcal{N})$ provide a $%
\mathbb{C}$-basis for the Hibi algebra $\mathcal{H}_{D(\mathcal{N})}$\ over $%
D(\mathcal{N})$ (\cite{Hi87}). The pattern corresponding to a $\mathcal{N}$%
-standard tableau consists of two parts having supports in $\{z_{b}^{(a)}\in
\digamma _{k,2n}:a\geq 2n\}$ and in $\{z_{b}^{(a)}\in \digamma _{k,2n}:a\leq
2n\}$ respectively, and therefore GT patterns of $GL_{k}$ and $GL_{2n}$ with
the same top row $D$. Furthermore, by Proposition \ref{support}, such GT
patterns of $GL_{2n}$ have their supports in $\{z_{b}^{(a)}\in \digamma
_{k,2n}:a\leq 2n$ and $b\leq (a+1)/2\}$, and then they represent GT\
patterns for $Sp_{2n}$ (e.g., \cite{Ki08}). Hence the correspondence given
in Proposition \ref{cone4GL} provides a bijection between the set of $%
\mathcal{N}$-standard monomials of shape $D$ and the set of pairs of
patterns $(\mathsf{p}_{-},\mathsf{p}_{+})$ where $\mathsf{p}_{-}$ and $%
\mathsf{p}_{+}$ are GT patterns for $GL_{k}$ and $Sp_{2n}$ respectively with
the same top row $D$. Therefore, we obtain the bijection between the set of
weight vectors for $\rho _{k}^{D}\otimes \sigma _{2n}^{D}$ and $\mathcal{P}(%
\mathcal{N}_{k,2n})_{D}$.
\end{proof}

\begin{example}
The GT pattern of $GL_{10}$ given in Example \ref{exampleSp}, considered as
an element of $\mathcal{P}_{4,6}=\mathcal{P}(4,6)/\left\langle \mathsf{p}_{%
\widehat{1}}\right\rangle $, can be visualized as a pattern over $\Gamma
_{4,6}$ as follows:%
\begin{equation*}
\begin{array}{ccccccccc}
&  &  & 5 &  &  &  &  &  \\ 
&  & 5 &  & 3 &  &  &  &  \\ 
& 5 &  & 4 &  & 1 &  &  &  \\ 
5 &  & 4 &  & 2 &  & 0 &  &  \\ 
& 5 &  & 4 &  & 1 &  & 0 &  \\ 
&  & 5 &  & 3 &  & 0 &  & 0 \\ 
&  &  & 4 &  & 2 &  & 0 &  \\ 
&  &  &  & 3 &  & 0 &  &  \\ 
&  &  &  &  & 2 &  &  & 
\end{array}%
\end{equation*}%
Note that the non-zero entries are corresponding to the poset $\digamma
_{4,6}$. As we discussed in Remark \ref{gluing} for $M_{n,m}$, it can also
be seen as the fiber product of two GT patterns, one for $GL_{4}$ and the
other for $Sp_{6}$, along their top rows:%
\begin{equation*}
\text{ }%
\begin{array}{ccccccc}
5 &  & 4 &  & 2 &  & 0 \\ 
& 5 &  & 4 &  & 1 &  \\ 
&  & 5 &  & 3 &  &  \\ 
&  &  & 5 &  &  & 
\end{array}%
\text{ and }%
\begin{array}{cccccc}
5 &  & 4 &  & 2 &  \\ 
& 5 &  & 4 &  & 1 \\ 
&  & 5 &  & 3 &  \\ 
&  &  & 4 &  & 2 \\ 
&  &  &  & 3 &  \\ 
&  &  &  &  & 2%
\end{array}%
\end{equation*}%
The above GT patterns correspond to the following semistandard tableaux of
shape $(5,4,2)$ via the conversion procedure given in Lemma \ref%
{pattern2tableau}:%
\begin{equation*}
\young(11111,2223,34)\text{ and }\young(11234,3345,56)
\end{equation*}%
They were denoted by $T^{-}$ and $T^{+}$ respectively in the proof of
Theorem \ref{standardmonomialN}, and they can be easily read from the row
indices and the column indices of the $\mathcal{N}$-standard monomial in
Example \ref{exampleSp}.
\end{example}

\medskip

Finally, we remark that the discussion in \cite{How05} on a simplicial
decomposition of a polyhedral cone and its relation with an algebra
decomposition can be directly applied to our case, and then we can interpret
standard monomial theory of $\mathcal{R}(\mathcal{N}_{k,2n})$ in terms of a
simplicial decomposition of $\mathcal{C}(\mathcal{N}_{k,2n})$. More
precisely, if we take a maximal linearly ordered subset $\mathcal{S}$ of $D(%
\mathcal{N})$, then all the products of elements from%
\begin{equation*}
\widehat{\mathcal{S}}=\{\delta _{\lbrack I:J]}\in \mathcal{R}(\mathcal{N}%
_{k,2n}):[I:J]\in \mathcal{S}\}
\end{equation*}%
are $\mathcal{N}$-standard monomials. Therefore, elements in $\widehat{%
\mathcal{S}}$ are algebraically independent and generate a polynomial
subring of $\mathcal{R}(\mathcal{N}_{k,2n})$. On the other hand, a
computation of Lemma \ref{pattern2tableau} shows that $\mathcal{S}$ must be
induced from a linearization of the poset $\digamma _{k,2n}$. Moreover, from
our construction of $\mathcal{C}(\mathcal{N}_{k,2n})$ in terms of $\digamma
_{k,2n}$, all possible linearizations of $\digamma _{k,2n}$ give rise to a
simplicial decomposition of the cone $\mathcal{C}(\mathcal{N}_{k,2n})$.
Consequently, we can obtain a decomposition of $\mathcal{R}(\mathcal{N}%
_{k,2n})$ into polynomial rings. This decomposition is not disjoint,
however, it is compatible with a simplicial decomposition of $\mathcal{C}(%
\mathcal{N}_{k,2n})$ induced from linearizations of $\digamma _{k,2n}$. For
more details in this direction, we refer readers to \cite{How05}.

\medskip

\noindent \textbf{Acknowledgment.} The author thanks Roger Howe for helpful
conversations on this project and Victor Protsak for insightful discussions
on the poset operations in Section \ref{CMnm}.

\end{document}